\documentclass[10pt]{article}

\usepackage{a4wide}
\usepackage{amssymb}
\usepackage{amsfonts}
\usepackage{amsmath}
\input xy
\xyoption{arrow} \xyoption{matrix}

\date{}

\newtheorem{proposition}{Proposition}[section]
\newtheorem{theorem}[proposition]{Theorem}
\newtheorem{lemma}[proposition]{Lemma}

\newtheorem{corollary}[proposition]{Corollary}

\def\der{\partial }

\def\nFM0{{\nu }_{F,M_0}}
\def\nFN0{{\nu }_{F,N_0}}
\def\nGN0{{\nu }_{G,N_0}}

\def\N0{ {\bf N}_0 }

\def\v{\varphi}
\def\ra{\rightarrow}

\def\Xpm{X^{\pm }}

\def\l1{{\lambda}_1}

\def\a{\alpha}
\def\a0{ {\alpha }_0}
\def\a1{ {\alpha }_1}

\def\l{\lambda}

\def\nFGM0{{\nu }_{F,G,M_0}}

\def\nFN0{{\nu}_{F,N_0}}

\def\sm{{\sigma}^m}

\def\sm1{{\sigma}^{-1}}

\def\smtp1{{\sigma}^{-t+1}}

\def\S1{S^{-1}}

\def\Xpm1{X^{\pm 1}_1}

\def\sPM1{{\sigma }^{\pm 1}}
\def\sMP1{{\sigma }^{\mp 1 }}

\def\di{{\rm d.ind}}

\def\L{\Lambda}

\def\Ytm1{Y^{t-1}}
\def\Yim1{Y^{i-1}}

\def\bQ{\overline{Q}}

\def\ker{ {\rm ker } }

\def\gr{ {\rm gr} }

\def\SL2Z{ {\rm SL}_2({\bf Z}) }

\def\th{ \theta }

\def\Gp1{ G^{1 , 1 } }
\def\P11{ P^{-1 , 1 } }
\def\Pp1{ P^{1 , 1 } }

\def\th{\theta}

\def\nCLsr{{}^\nu\kern-2pt {\cal L}^{\sigma , \rho  }}
\def\nP{{}^\nu \kern-2pt P}
\def\nL{{}^\nu\kern-2pt L}
\def\nLL{{}^\nu\kern-2pt \Lambda}
\def\nPsr{{}^\nu\kern-2pt P^{\sigma , \rho  }}
\def\nLsr{{}^\nu\kern-2pt L^{\sigma , \rho  }}
\def\nuCL{{}^\nu\kern-2pt  {\cal L}}
\def\nCLsr{{}^\nu\kern-2pt {\cal L}^{\sigma , \rho  }}
\def\nCL1m{{}^\nu\kern-2pt {\cal L}^{-1 , 1  }}
\def\x1nu{x^\frac{1}{\nu}}
\def\xm1nu{x^{-\frac{1}{\nu}}}

\def\ra{\rightarrow }

\def\CB{{\cal B}}

\def\CC{ {\cal C}}

\def\nAM0{{\nu }_{{\cal A},M_0}}
\def\nAN0{{\nu }_{{\cal A},N_0}}

\def\bR{\overline{R}}

\def\bQ{\overline{Q}}

\def\bx{\overline{x}}

\def\ga{\mathfrak{a}}
\def\gb{\mathfrak{b}}

\def\gn{\mathfrak{n}}
\def\gm{\mathfrak{m}}
\def\gp{\mathfrak{p}}

\def\gr{\mathfrak{r}}

\def\SL{{\rm SL}}

\def\di!{\frac{\der^i}{i!}}
\def\dik!{\frac{\der^k_i}{k!}}

\def\Max{{\rm Max}}

\def\N{\mathbb{N}}

\def\0{\overline{0}}
\def\1{\overline{1}}

\def\Ln1{\L_{n,\overline{1}}}

\def\oa{\overline{a}}

\def\a1{a_{\overline{1}}}

\def\bs{\overline{s}}

\def\S{\Sigma}

\def\vn1{\overrightarrow{n-1}}

\def\gf{\mathfrak{f}}

\def\Min{{\rm Min}}

\def\bu{\overline{u}}

\def\mJ{\mathbb{J}}
\def\mI{\mathbb{I}}

\def\K1{{\rm K}_1}

\def\hmI1{\widehat{\mI_1}}
\def\tmI1{\widetilde{\mI_1}}
\def\tmJ1{\widetilde{\mJ_1}}
\def\hB1{\widehat{B_1}}
\def\hCB1{\widehat{\CB_1}}

\def\OCC{\overline{\CC}}
\def\br{\overline{r}}
\def\bc{\overline{c}}
\def\bn{\overline{n}}
\def\gt{\mathfrak{t}}
\def\bs{\overline{s}}
\def\bt{\overline{t}}
\def\ga{\mathfrak{a}}

\def\tor{{\rm tor}}
\def\udim{{\rm udim}}

\begin{document}

\author{V. V. \  Bavula 
}

\title{Characterizations of left orders in left Artinian rings}

\maketitle

\begin{abstract}
Small (1966),  Robson (1967), Tachikawa (1971) and Hajarnavis (1972)  have given different criteria for a ring to have a left Artinian left quotient ring.  In this paper, three   more new criteria are given.




$${\bf Contents}$$
\begin{enumerate}
\item Introduction.
\item Necessary and sufficient conditions for a ring to have a left Artinian left quotient ring.
\item Associated graded ring.
 \item Criteria similar to  Robson's Criterion.
 \item A left quotient ring of a factor ring.
\end{enumerate}
\end{abstract}


\section{Introduction}

In this paper, module means a left module, and the following notation is fixed:
\begin{itemize}
\item  $R$ is a ring with 1;
\item   $\CC = \CC_R$  is the set of regular elements of the ring $R$ (i.e. $\CC$ is the set of non-zero-divisors of the ring $R$);
\item   $Q=Q_{l,cl}(R):= \CC^{-1}R$ is the {\em left quotient ring}  (the {\em classical left ring of fractions}) of the ring $R$ (if it exists);
\item   $\gn$ is a prime radical of $R$ and $\nu$ is its {\em nilpotentcy degree} ($\gn^\nu \neq 0$ but $\gn^{\nu +1}=0$);
\item   $\bR := R/ \gn$ and $\pi: R\ra \bR$, $r\mapsto \br =r+\gn$;
\item   $\OCC := \CC_{\bR}$ is the set of regular elements of the ring $\bR$ and $\bQ := \OCC^{-1}\bR$ is its left quotient ring;
\item   $\CC':= \pi^{-1}(\OCC):=\{ c\in R\, \, | \, c+\gn \in \OCC\}$ and $Q':=\CC'^{-1}R$ (if it exists).
\end{itemize}
Goldie's Theorem \cite{Goldie-PLMS-1960} characterizes left orders in semi-simple rings, it is a criterion of when the left quotient ring of  a ring is semi-simple (earlier,  characterizations were given, by Goldie \cite{Goldie-PLMS-1958} and Lesieur and Croisot \cite{Lesieur-Croisot-1959}, of left orders in a simple Artinian ring). Talintyre \cite{Talintyre-QuotRinMax-63} and Feller and Swokowski \cite{Feller-Swokowski-RefPROCams-61} have given conditions
which are sufficient for a left Noetherian ring to
have a left quotient ring. Further, for a left Noetherian ring which has
a left quotient ring, Talintyre \cite{Talintyre-QuotRinMin-66} has established necessary and sufficient
conditions for the left quotient ring to be left Artinian. Small \cite{Small-ArtQuotRings-66, Small-CorrArtQuotRings-66},    Robson \cite{Robson-ArtQuotRings-67}, and latter  Tachikawa \cite{Tachikawa-AutQuotR-71} and Hajarnavis \cite{Hajarnavis-ThmSmall-72}  have given different criteria  for a ring to have a left Artinian left quotient ring.  In this paper, three   more new criteria are given (Theorem \ref{2Sep12}, Theorem \ref{3Sep12} and Theorem \ref{4Sep12}). Theorem \ref{13Sep12} gives an affirmative  answer to the question: {\em Let $R$ be a ring with a left Artinian left quotient ring $Q$ and $I$ be a $\CC$-closed ideal of $R$ such that $I\subseteq \gn$. Is the left quotient ring $Q(R/I)$ of $R/I$ a left Artinian ring?}

In the proofs of all the criteria (old and new) Goldie's Theorem is used. Each of the criteria comprises  several conditions. The conditions in the criteria of Small, Robson and Hajarnavis are `strong' and are given in terms of the ring $R$ rather than of its factor $\bR = R/ \gn $. On the contrary,  the conditions of the criteria of the present paper are `weak' and given in terms of the ring $\bR$ and the {\em finite} set of explicit $\bR$-modules (they are subfactors of the prime radical $\gn$ of the ring $R$). As a result the proofs in the present paper are different from the proofs of the earlier criteria and different ideas are used. since the conditions of the criteria of the present paper are weak.  One can produce more criteria that are `intermediate' between these and  the older ones (and several such `derivatives' are given). There is a logic and practical need to have each type of criteria -- `weak'  and `strong'. Weak ones are used to check that a ring has a left Artinian left quotient ring. Then the strong ones are used to produce strong corollaries about its properties.

The set $\CC$ is not always a left (or right) Ore set in $R$ as a result the classical left quotient ring $\CC^{-1}R$ does not always exist but there {\em always} exists the {\em largest left Ore set} $\CC_{l,R}$ and the {\em largest right Ore set} $\CC_{r, R}$ in $\CC$ of the ring $R$, \cite{larglquot}. In general, $\CC_{l,R}\neq \CC_{r, R}$, \cite{Bav-intdifline}. In \cite{larglquot}, the {\em largest left quotient ring } $Q_l(R):= \CC_{l,R}^{-1}R$ and the {\em largest right quotient ring} $Q_r(R):=R \CC_{r,R}^{-1}$ are introduced. In  \cite{Bav-intdifline}, these rings are found for the ring $\mI_1= K\langle x, \der , \int\rangle $ of polynomial integro-differential operators over a field $K$ of characteristic zero, $\CC_{l,\mI_1}\neq \CC_{r,\mI_1}$ and $Q_l(R)\not\simeq  Q_r(R)$.



\section{Necessary and sufficient conditions for a ring to have a left Artinian left quotient ring}\label{NCLALQ}

The aim of this section is to prove Theorem \ref{2Sep12} which is a criterion for a ring $R$ to have   a left Artinian left quotient ring. We give a Recipe to produce `derivative' criteria (see Corollary \ref{a9Sep12} and Theorem \ref{C2Sep12}). Using Theorem \ref{2Sep12} and Theorem \ref{C2Sep12} in combination with results of Small and P. F. Smith, criteria are obtained for a left Noetherian ring $R$ (Corollary \ref{A2Sep12}) and for a commutative ring $R$ (Corollary \ref{x10Sep12}) to have left Artinian left quotient ring.

Let us recall   certain properties of left Artinian rings.
\begin{proposition}\label{}
{\rm (Proposition 3.1, \cite{Aus-Rei-Smalo-book-1995})}  Let $A$ be a left Artinian ring and $\gr$ be its radical. Then
\begin{enumerate}
\item The radical $\gr$ of $A$ is a nilpotent ideal.
\item The factor ring $A/ \gr$ is a semi-simple.
\item An $A$-module $M$ is semi-simple iff $\gr M=0$.
\item There is only finite number of non-isomorphic simple $A$-modules.
\item The ring $A$ is a left Noetherian ring.
\end{enumerate}
\end{proposition}

A ring $R$ is called a {\em left Goldie ring} if it satisfies ACC (the {\em ascending chain condition}) for left annihilators  and contains no infinite direct sums of left ideals.

Suppose that a ring $R$ satisfies the condition (a) of Theorem \ref{2Sep12}, i.e. $\bR$ is a {\em (semiprime) left Goldie ring}. By Goldie's Theorem, its  left quotient ring $\bQ := \OCC^{-1}\bR$ is a semisimple (Artinian) ring where $\OCC$ is the set of  regular elements of the ring $\bR$. The ring $R$ admits {\em the} $\gn$-{\em adic filtration} ({\em the prime radical filtration}):
\begin{equation}\label{chid}
\gn^0:=R\supset\gn \supset \cdots \supset \gn^i\supset \cdots
\end{equation}
which stops at $(\nu +1)$'th step if $\gn^\nu\neq 0$ but $\gn^{\nu +1}=0$, i.e. $\nu$ is the nilpotency degree of the ideal $\gn$. The {\em associated graded algebra}
 $${\rm gr} \, R= \bR\oplus\gn / \gn^2\oplus \cdots \gn^i/ \gn^{i+1}\oplus \cdots $$ is an $\N$-graded ring and every component $\gn^i/\gn^{i+1}$ is an $\bR$-bimodule. Recall that $\OCC$ is a left Ore set in $\bR$ (by Goldie's Theorem) and that module means a left module. For each integer $i\geq 1$, let
\begin{equation}\label{ttidef}
\tau_i:=\tor_{\OCC}(\gn^i/\gn^{i+1}):=\{ u\in\gn^i/\gn^{i+1}\, | \, \bc u =0\; \; {\rm for \; some}\; \bc \in \OCC\}
\end{equation}
be the $\OCC$-{\em torsion submodule} of the {\em left} $\bR$-module $\gn^i/\gn^{i+1}$. Clearly, $\tau_i$ is an $\bR$-bimodule. Then the $\bR$-bimodule
\begin{equation}\label{ttidef1}
\gf_i:=(\gn^i/\gn^{i+1})/\tau_i
\end{equation}
is a $\OCC$-{\em torsionfree} left $\bR$-module. There is a unique  ideal, say $\gt_i$, of the ring $R$ such that
$$\gn^{i+1}\subseteq \gt_i\subseteq \gn^i\;\; {\rm and}\;\; \gt_i/ \gn^{i+1}=\tau_i.$$ Clearly, $\gf_i\simeq \gn^i/ \gt_i$.

\begin{theorem}\label{2Sep12}
Let $R$ be a ring. The following statements are equivalent.
\begin{enumerate}
\item The ring $R$ has  a left Artinian left quotient ring $Q$.
\item
\begin{enumerate}
\item The ring $\bR$ is a left Goldie ring $Q$.
\item $\gn$ is a nilpotent ideal.
\item $\CC'\subseteq \CC$.
\item  The left $\bR$-modules $\gf_i$, $1\leq i\leq 1$, contain no infinite direct sums of nonzero submodules, and
\item for every element $\bc\in \OCC$, the map $\cdot \bc: \gf_i\ra \gf_i$, $f\mapsto f\bc$, is an injection; equivalently, if, for an element $a\in \gn^i/ \gn^{i+1}$, there are elements $\bs$, $\bc\in \OCC
    $ such that $\bs a \bc =0$ then $\bt a=0$ for some element $\bt\in \OCC$;  equivalently, if, for an element $a\in \gn^i/ \gn^{i+1}$, there is an  element $\bc\in \OCC $ such that $ a \bc =0$ then $\bt a=0$ for some element $\bt\in \OCC$.
\end{enumerate}
\end{enumerate}
If one of the equivalent conditions holds then $\CC = \CC'$, $\CC^{-1}\gn$ is the prime radical of the ring $Q$ which is a nilpotent ideal of nilpotency degree $\nu$, and the map $Q/ \CC^{-1}\gn\ra \bQ$, $c^{-1}r\mapsto \bc^{-1}\br$, is a ring isomorphism with the inverse $\bc^{-1}\br\mapsto c^{-1}r$.
\end{theorem}

{\it Proof}. $(2\Rightarrow 1)$  The proof comprises three steps:

(A) $\CC'$ {\em is a left Ore set in} $R$. From which it follows easily that

(B) {\em the ring $Q':= \CC'^{-1}R$ is a left Artinian ring}, and then that

(C) $\CC= \CC'$ {\em and} $Q=Q'$.

The key point of the proof of step (A) is the following well known  lemma ($\udim$ is the {\em uniform dimension}; for the properties of the uniform dimension the interested reader is referred to \cite{Jategaonkar-LocNRings} or \cite{MR}).

\begin{lemma}\label{a2Sep12}
Let $\bR$ be a left Goldie ring and $M\subseteq N$ be $\OCC$-torsionfree $\bR$-modules such that  $\udim_{\bR}(M) = \udim_{\bR}(N)<\infty$. Then the $\bR$-module $N/M$ is $\OCC$-torsion.
\end{lemma}

(A) $\CC'$ {\em is a left Ore set in} $R$. We have to show that, for elements $c\in \CC'$ and $r\in R$, there exist elements $c'\in \CC'$ and $r'\in R$ such that
$$ c'r=r'c.$$
To prove this statement we use a downward induction on the {\em degree} $i= \deg (r)$ of the element $r$ which is, by definition,  the unique number $i\in \{ 0,1, \ldots , \nu +1\}$ such that $r\in \gn^i\backslash \gn^{i+1}$ if $r\neq 0$ and $\deg (0):=\nu +1$. The base of the induction $i=\nu +1$, i.e. $r=0$, is obvious (take $c'=1$ and $ r'=0$). Suppose that $i<\nu +1$ and that the statement (A) is true for all elements of the ring $R$ with degrees $i'$ such that $i<i'\leq \nu +1$.

{\em Case 1:} $i=0$. The ring $\bR$ is a semiprime left Goldie ring. By Goldie's Theorem,
$$ \bc_1\br =\br_1\bc$$ for some elements $\bc_1= c_1+\gn \in \OCC$ and $\br_1=r_1+\gn\in \bR$ where, for $x\in \bR$, $\bx := x+\gn$. By the condition (c), $c_1\in \CC'$. Then, the element
$$ n:= c_1r-r_1c\in \gn$$
has degree $\geq 1$, and, by induction, there are elements $c_2\in \CC'$ and $r_2\in R$ such that $c_2n=r_2c$, and so
$$ c_2c_1r=(r_2+c_2r_1)c.$$
It suffices to take $c'=c_2c_1$ and $r'=r_2+c_2r_1$.

{\em Case 2:} $i>0$. The left $\bR$-module $\gf_i$ is a $\OCC$-{\em torsionfree} $\bR$-module with $\udim_{\bR}(\gf_i)<\infty$ (by the condition (d)),  and $r+\tau_i\in \gf_i$. Then, by the condition (e), the map $\cdot c:\gf_i\ra \gf_i$, $f\mapsto fc$, is an $\bR$-module monomorphism since $ fc= f\bc$ where $f\in \gf_i$ and $\bc = c+\gn \in \OCC$. In particular,
$$ \udim_{\bR}(\gf_ic) = \udim_{\bR}(\gf_i)<\infty.$$
By Lemma \ref{a2Sep12}, the $\bR$-module $\gf_i / \gf_ic$ is $\OCC$-{\em torsion}. In particular, by the condition (c), there exist element $s\in \CC'$ and $a\in R$ such that
$$ (sr-ac)+\gn^{i+1}\in \tau_i.$$
Fix an element $t\in \CC'$ such that $\l  := t(sr-ac)\in \gn^{i+1}$, and so $\deg (\l ) \geq i+1$. By induction, there are elements $l\in \CC'$ and $b\in R$ such that $l\l = bc$, and so
$$ ltsr=(lta+b) c.$$
It suffices to take $c'=lts$ and $r'=lta+b$. The proof of statement (A) is complete.

Let us deduce obvious corollaries of (A) (where $\gn\CC'^{-1}:=\{ nc'^{-1}\, | \, n\in \gn , c'\in \CC'\}$):

(A1) $\gn\CC'^{-1}\subseteq \CC'^{-1}\gn$. In particular, $\CC'^{-1}\gn$ is an {\em ideal}  of the ring $Q'$ such that

(A2) $(\CC'^{-1}\gn )^i=\CC'^{-1}\gn^i$  {\em for all } $i\geq 1$ (obvious, see (A1)), and

(A3) {\em the map}
$$\CC'^{-1}R/\CC'^{-1}\gn\ra \bQ, \;\; c^{-1}r\mapsto \bc^{-1} \br, $$ {\em is a ring isomorphism with the inverse} $\bc^{-1} \br\mapsto c^{-1}r$.

{\it Proof of (A3)}: By (A1), the condition (c) and the universal property of left Ore localization, the both
 ring homomorphisms in (A3) are well-defined; clearly they are mutually inverse.

 (A4) $\CC'^{-1}\gn$ {\em is the prime radical of the ring} $Q'$ (obvious, see (A1) and (A3)).

  {\it Proof of (A1)}: We have to show that, for elements $n\in \gn$ and $c\in \CC'$, there exist elements $ n_1\in \gn$ and $ c_1\in \CC'$ such that $$ nc^{-1} = c_1^{-1} n_1.$$
 The equality is equivalent to the equality $c_1n=n_1c$. By (A), such an equality exists but with $n_1\in R$ rather than with $n_1\in \gn$. It suffices  to show that $n_1\in \gn$.  Since $0=\bc_1 \bn = \bn_1\bc$ in $\bR$, we conclude that  $\bn_1=0$ since $\bc\in \OCC$, i.e. $n_1\in \gn$. The proof of the statement (A1) is complete.

(B) {\em The ring $Q'$ is a left Artinian ring}.

Localizing at $\CC'$ the chain of ideals (\ref{chid}) of the ring $R$  we obtain the chain of ideals (by (A2)) of the ring $Q'$:
\begin{equation}\label{chid1}
\CC'^{-1} \gn^0=Q'\supset\CC'^{-1}\gn \supset \cdots \supset \CC'^{-1}\gn^i\supset \cdots \supseteq \CC'^{-1}\gn^\nu \supseteq \CC'^{-1}\gn^{\nu +1}=0.
\end{equation}
By (A1)--(A3), the factors $$\CC'^{-1}\gn^i/\CC'^{-1}\gn^{i+1}=(\CC'^{-1}\gn )^i/(\CC'^{-1}\gn )^{i+1}$$ are $Q'$-modules and $\bQ$-modules. By (A3), these modules structures are the `same' in the sense that the latter is induced by the former.  The statement (B) follows from the fact that all the factors $\CC'^{-1}\gn^i/\CC'^{-1}\gn^{i+1}$ are left $Q'$-modules/$\bQ$-modules of finite length. In more detail, the case $i=0$ is obvious (see (A3) and the condition (a)). For all $i=1, \ldots , \nu$,
$$ \CC'^{-1}\gn^i/\CC'^{-1}\gn^{i+1}\simeq \CC'^{-1}(\gn^i/\gn^{i+1})\simeq
\CC'^{-1}\gf^i\simeq \OCC^{-1}\gf_i,$$
the last isomorphism is due to (A3). By the statements (a) and  (d), the $\bQ$-module $\OCC^{-1}\gf_i$ is semisimple with
$$ \udim_{\bQ} (\OCC^{-1}\gf_i) = \udim_{\bR}(\gf_i)<\infty.$$
This means that the $\bQ$-modules $\OCC^{-1}\gf_i$ have finite length. This finishes the proof of statement (B).

Let $Q'^*$ be the group of units of the ring $Q'$. Let us prove one easy corollary of the statement (B) from which the statement (C) will follows easily.

(B1): $\CC \subseteq Q'^*$.

Recall that any monomorphism of a {\em finite length} module  is an automorphism. By (B), for each element $c\in \CC$, the $Q'$-module homomorphism $\cdot c: Q'\ra Q'$, $q\mapsto qc$, is a monomorphism, hence an automorphism necessarily of the type $\cdot d: Q'\ra Q'$, $q\mapsto qd$ for some element $d\in Q'$. Then $d=c^{-1}$, i.e. $\CC \subseteq Q'^*$.

(C) $\CC= \CC'$ {\em and} $Q=Q'$.

It suffices to establish only the first equality as it implies the second. Let $c\in \CC$. In view of the condition (c), it suffices  to show that $\bc = c+\gn \in \OCC$ (since then $c\in \CC'$). By (B1), $c\in Q'^*$, and, then, by (A3), the element $c+\CC'^{-1}\gn$ of the ring $ Q'/\CC'^{-1} \gn\simeq \bQ$ is a unit. Equivalently, the element $c+\gn\in \bR$, which is the image of the element $c+\CC'^{-1}\gn$ under the isomorphism in (A3),  is a unit of the ring $\bQ$. In particular, $c+\gn \in \OCC$, as required.  Therefore, the proof of the implication $(1\Rightarrow 2)$ is complete. Summarizing, by (A1)-(A4) and (C),  $\CC = \CC'$, $\CC^{-1}\gn$ is the prime radical of the ring which is a nilpotent radical of nilpotency degree $\nu$, and the map $Q/ \CC^{-1}\gn\ra \bQ$, $c^{-1}r\mapsto \bc^{-1}\br$, is a ring isomorphism with the inverse $\bc^{-1}\br\mapsto c^{-1}r$.

 $(1\Rightarrow 2)$  See Corollary \ref{a14Sep12}, this is a direct proof of the implication  $1\Rightarrow 2$.

 Alternatively, the conditions (a) and (b) (resp. (c)) are present in Robson's Criterion,   Theorem \ref{Th-2.10-Robson} (resp. Small's Criterion, Theorem C, \cite{Small-CorrArtQuotRings-66}). The condition that `the ring $R$ is $\gn$-quorite' in Robson's Criterion implies the condition (e). The conditions in Robson's Criterion  that `$R$ is $\gn$-reflective' (i.e. $\CC' = \CC$) and `$R$ satisfies ACC for $\CC$-closed ideals' imply the condition (d).   $\Box $


\begin{corollary}\label{CA2Sep12}
Let $R$ be a commutative Noetherian ring. The following  statements are equivalent.
\begin{enumerate}
\item The ring $R$ has an Artinian  quotient ring.
\item The conditions (c) of Theorem \ref{2Sep12} holds.
\item The associated primes of $(0)$ are the minimal primes of the ring $R$.
\end{enumerate}
\end{corollary}
{\it Proof}. $(1\Leftrightarrow 2)$ Theorem \ref{2Sep12}.

 $(2\Leftrightarrow 3)$ This equivalence was established by Small (Theorem 2.13, \cite{Small-ArtQuotRings-66} and Theorem C, \cite{Small-CorrArtQuotRings-66}).

  $(1\Leftrightarrow 3)$ Robson (Theorem 2.11, \cite{Robson-ArtQuotRings-67}). $\Box $

$\noindent $

{\bf Recipe for producing new criteria for a ring to have a left Artinian left quotient ring}: Take a criterion, for example Theorem \ref{2Sep12}, and replace any set $S$ of conditions by a {\em stronger} set of conditions $S'$ (i.e. $S'\Rightarrow S$) that are satisfied by the rings having left Artinian left quotient ring, and we obviously have a new criterion (the conditions of the new criterion are obviously necessary and, since they imply the old criterion, sufficient).  The weaker the original conditions, like in Theorem \ref{2Sep12}, the more we can produce. Precisely, using this Recipe, Corollary \ref{a9Sep12} and  Corollary \ref{A3Sep12} are proved. The conditions of Robson and Small are `strong'. So, it seems that the Recipe does not work in these cases but we can go in the opposite direction -- to weaken some of the original conditions but usually it is a hard work to prove things in this case. Precisely on this way a `weaker' version of Robson's Criterion is obtained (Theorem \ref{4Sep12}).

Let $M$ be a left $\bR$-module. A submodule $N$ of $M$ is called a $\OCC$-{\em closed} if, for elements $m\in M$ and $\bc\in \OCC$, $\bc m\in N$ implies $m\in N$ (equivalently, $\tor_{\OCC}(M/N)=0$).

\begin{corollary}\label{a9Sep12}
Let $R$ be a  ring. The following two statements are equivalent.
\begin{enumerate}
\item The  ring $R$ has  a left Artinian left quotient ring.
\item The conditions (a)--(c) and (e) of Theorem \ref{2Sep12}  and any of the following conditions hold,

    (d1) the left $\bR$-modules $\gf_i$, where $i\geq 1$, satisfy ACC on $\OCC$-closed submodules,

(d2) the left $\bR$-modules $\gf_i$, where $i\geq 1$, satisfy DCC on $\OCC$-closed submodules.
\end{enumerate}
\end{corollary}

{\it Proof}. $(1\Rightarrow 2)$  If statement 1 holds then every left $\bQ$-module $\OCC^{-1}\gf_i$ is semisimple of finite length. The $\bR$-module $\gf_i$ is $\OCC$-torsionfree. Therefore,  $\gf_i\subseteq \OCC^{-1}\gf_i$ and the map $I\mapsto \OCC^{-1}I$ is a bijection between the set of $\OCC$-closed submodules of $\gf_i$ and the $\bQ$-submodules of $\OCC^{-1}\gf_i$ with the inverse map $J\mapsto \gf_i\cap J$. Now, it is obvious that the conditions (d1) and (d2) are satisfied.

 $(2\Rightarrow 1)$ In view of the bijection $I\mapsto \OCC^{-1}I$, each of the conditions (d1) or (d2) implies the condition (d) of Theorem \ref{2Sep12}. $\Box$

 $\noindent $

A ring $R$ is called $\gn$-{\em reflective} if, for $c\in R$, $c\in \CC$ iff $c+\gn \in \OCC$; equivalently, $\CC'=\CC$.

\begin{corollary}\label{A3Sep12}
Let $R$ be a left Noetherian ring. The following two statements are equivalent.
\begin{enumerate}
\item The  ring $R$ has  a left Artinian left quotient ring.
\item The ring $R$ is $\gn$-reflective and the conditions (a), (b), (d) and (e) of Theorem \ref{2Sep12} hold.
\end{enumerate}
\end{corollary}

{\it Proof}. The condition `$R$ is $\gn$-reflective', i.e. $\CC'=\CC$, is present in Robson's Criterion (Theorem 2.10, \cite{Robson-ArtQuotRings-67}) and it is stronger that the condition (c) of Theorem \ref{2Sep12}. So, replacing the latter by the former in Theorem \ref{2Sep12} and using the Recipe we finish the proof of the corollary.  $\Box $

$\noindent $

For a ring $R$ having a left Artinian left quotient ring $Q$, Theorem \ref{C2Sep12} provides many examples of left Ore subsets $\CC''\subseteq \CC$ such that $\CC''^{-1}R\simeq Q$.

\begin{theorem}\label{C2Sep12}
Let $R$ be a ring. The following statements are equivalent.
\begin{enumerate}
\item The  ring $R$ has  a left Artinian left quotient ring $Q$.
\item The conditions (a), (b), (c'), (d) and (e) hold  (see Theorem \ref{2Sep12}) where

(c') for each element $\alpha \in \OCC$, there exists a regular element $c=c(\alpha )\in \CC$ such that $\alpha = c+\gn$; equivalently,  there exists a submonoid $\CC''\subseteq \CC$ such that $\overline{\CC''}= \OCC$ (where $\overline{\CC''}:=\{ c+\gn \, | \, c\in \CC''\})$.
\end{enumerate}
If one of the equivalent conditions holds then $\CC''$ is a left Ore set in $R$, $\CC''^{-1}R=Q$,  $\CC''^{-1}\gn$ is the prime radical of the ring $Q$ which is a nilpotent ideal of nilpotency degree $\nu$, and the map $Q/ \CC''^{-1}\gn\ra \bQ$, $c^{-1}r\mapsto \bc^{-1}\br$, is a ring isomorphism with the inverse $\bc^{-1}\br\mapsto c^{-1}r$ where $c$ is any element of $\CC''$ such that $ \bc = c+\gn$.
\end{theorem}

{\it Proof}. $(1\Rightarrow 2)$ Take $\CC''= \CC'$. Then the implication $(1\Rightarrow 2)$ follows from Theorem \ref{2Sep12}.

$(2\Rightarrow 1)$  The proof is almost identical to the proof of the implication $(2\Rightarrow 1)$ in Theorem \ref{2Sep12}. Statements (A), (B), (A1)--(A4) and (B1) hold where $\CC'$ and $Q'$ are replaced by $\CC''$ and $Q'':= \CC''^{-1}R$ respectively. Let us write down these statements for future references  adding the superscript $()''$.

(A'') $\CC''$ {\em is a left Ore set in} $R$.

(B'') {\em the ring $Q''$ is a left Artinian}.

(A1'') $\gn\CC''^{-1}\subseteq \CC''^{-1}\gn$.

(A2'') $(\CC''^{-1}\gn )^i=\CC''^{-1}\gn^i$  {\em for all } $i\geq 1$.

(A3'') {\em The map}
$$\CC''^{-1}R/\CC''^{-1}\gn\ra \bQ, \;\; c^{-1}r\mapsto \bc^{-1} \br, $$ {\em is a ring isomorphism with the inverse} $\bc^{-1} \br\mapsto c^{-1}r$ {\em where $c$ is any element of $\CC''$ such that} $\bc = c+\gn$.

(A4'') $\CC''^{-1}\gn$ {\em is the prime radical of the ring} $Q''$.

(B1''): $\CC \subseteq Q''^*$ {\em where} $Q''^*$ {\em is the group of units of the ring} $Q''$.

The proofs of the above statements are identical (with trivial obvious modifications) to the proofs of their counterparts and are left to the reader  as an easy exercise. The statement (B1'') admits a more direct proof:  Let $c\in \CC$. By the condition (c'), $c=s+n$ for some elements $s\in \CC''$ and $n\in \gn$. Then $c=(1+ns^{-1})s= (1-t^{-1}m)s$ for some elements $t\in \CC''$ and $n\in \gn$ (by (A1''), $ns^{-1} = -t^{-1}m$). By (A2''), the element $t^{-1}m$ is a nilpotent element, i.e. $(t^{-1}m)^{k+1}=0$ for some $k\geq 0$. Then
$$ c^{-1} = s^{-1}(1-t^{-1}m)^{-1}= s^{-1}(1+\sum_{i=1}^k(t^{-1}m)^i)= s^{-1}(1+\sum_{i=1}^kt_i^{-1}m_i)\in Q'',$$
 for some elements $t_i\in \CC''$ and $m_i\in R$.

(B2'') $\CC$ {\em is a left Ore set in} $R$. For elements $c\in \CC$ and $r\in R$, we have to find elements $c_1\in \CC$ and $r_1\in R$ such that $c_1r=r_1c$. By (B1''), $c\in Q''^*$ and so $Q''\ni rc^{-1}= c_1^{-1}r_1$ for some $c_1\in \CC''\subseteq \CC$ and $r_1\in R$; equivalently, $c_1r=r_1c$.


By (B2'') and (B1''), we conclude that

 (C'') $Q=Q''$.

 In particular, the ring $Q$ is a left Artinian ring, by (B'') and (C'').  $\Box$

\begin{corollary}\label{A2Sep12}
Let $R$ be a left Noetherian ring. The following two statements are equivalent.
\begin{enumerate}
\item The ring $R$ has a left Artinian left quotient ring.
\item $\CC'\subseteq \CC $.
\item For each element $\alpha \in \OCC$, there exists an element $c=c(\alpha )\in \CC$ such that $\alpha = c+\gn$.
\end{enumerate}
\end{corollary}

{\it Proof}. $(1\Leftrightarrow 2)$ This is due to Small \cite{Small-ArtQuotRings-66}.

$(2\Rightarrow 3)$ Trivial.

$(3\Rightarrow 1)$ It suffices to show that the conditions of Theorem \ref{C2Sep12} hold. The ring $R$ is left Noetherian, so the conditions (a), (b) and (d) hold.  The condition (c') is given. So, it remains to show that the condition (e) is true. Let $\CC''$ be the submonoid of $\CC$ generated by all the elements $\{ c(\alpha ) \, | \, \alpha \in \OCC \}$. By (Proposition 4.3, p. 288, \cite{Stenstrom-RingQuot}), $\CC''$ is a left Ore set. Then

(A1'') $\gn\CC''^{-1}\subseteq \CC''^{-1}\gn$.

Let $n\in \gn$ and $c''\in \CC'$. Then $ nc''^{-1}=c_1''^{-1} n_1$ for some elements $c_1''\in \CC''$ and $n_1\in R$ (since $\CC''$ is a left Ore set) or, equivalently, $c_1'' n=n_1c''$. We have to show that $n_1\in \gn$. Then $\bn_1 \overline{c''}=0$ in $\bR$, and so $\bn_1=0$ since $\overline{c''}\in \OCC$, i.e. $n_1\in \gn$.

(A2'') $\gn^i\CC''^{-1}\subseteq \CC''^{-1}\gn^i$, $i\geq 1$ (by (A1'')).

Then the statement (A2'') implies the statement (e): if $\oa \bc =0$ where $\oa = a+\gn^{i+1} \in \gn^i/\gn^{i+1}$ and $\bc = c+\gn \in \OCC$, i.e. $b:= ac\in \gn^{i+1}$ where $c\in \CC''$ (by statement 3), then, by (A2''), $bc^{-1} = c_1^{-1} b_1$ for some elements $c_1\in \CC''$ and $b_1\in \gn^{i+1}$, and so
$$c_1a=c_1(bc^{-1}) = c_1(c_1^{-1} b_1) = b_1\in \gn^{i+1}, $$
i.e. $\bc_1 \oa =0$ where $\bc_1\in \OCC$. $\Box $

\begin{corollary}\label{x10Sep12}
Let $R$ be a commutative  ring. The following  statements are equivalent.
\begin{enumerate}
\item The ring $R$ has an Artinian  quotient ring.
\item
\begin{enumerate}
\item The ring $\bR$ is a  Goldie ring.
\item $\gn$ is a nilpotent ideal.
\item $\CC'\subseteq \CC$.
\item  The  $\bR$-modules $\gf_i$, $1\leq i\leq \nu $, contain no infinite direct sums of nonzero submodules.
\end{enumerate}
\item
\begin{enumerate}
\item The ring $\bR$ is a  Goldie ring.
\item $\gn$ is a nilpotent ideal.
\item For each element $\alpha \in \OCC$, there exists an element $c=c(\alpha )\in \CC$ such that $\alpha = c+\gn$.
\item  The  $\bR$-modules $\gf_i$,  $1\leq i\leq \nu$, contain no infinite direct sums of nonzero submodules.
\end{enumerate}
\item $R$ is a Goldie ring and $\CC'\subseteq \CC $.
\item $R$ is a Goldie ring and, for each element $\alpha \in \OCC$, there exists an element $c=c(\alpha )\in \CC$ such that $\alpha = c+\gn$.
\end{enumerate}
\end{corollary}

{\it Proof}. $(1\Leftrightarrow 2)$ Theorem \ref{2Sep12}.

$(1\Leftrightarrow 3)$ Theorem \ref{C2Sep12}.

$(1\Leftrightarrow 4)$ This is due to P. F. Smith (Theorem 2.11, \cite{Hajarnavis-ThmSmall-72}).

$(4\Rightarrow 5)$ Trivial.

$(5\Rightarrow 4)$  The condition $\CC'\subseteq \CC$ is equivalent to two conditions: $\pi ( \CC ) = \OCC$ and $\CC +\gn \subseteq \CC$ where $\pi : R\ra \bR$, $r\mapsto \br$. By statement 5, the first condition is given. Let $c\in \CC$ and $n\in \gn$. To prove the second statement we have to show that $c+n\in \CC$. Notice that  $n$ is a nilpotent element and the ring $R$ is a subring of $\CC^{-1}R$. Now, the element $c+n= c(1+c^{-1}n)$ is a unit of the ring $\CC^{-1}R$ (as a product of two units). Therefore, $c+n\in \CC$.  $\Box$


\section{Associated graded ring}\label{ASGR}

The aim of this section is to give  another criterion (Theorem \ref{3Sep12}) for a ring $R$ to have a left Artinian left quotient ring via its associated graded ring ${\rm gr}\, R$ with respect to the $\gn$-adic  filtration.

A multiplicative  set $S$ of a ring $R$ is a {\em left denominator set} if it is a left Ore set and if $rs=0$, for some elements $r\in R$ and $s\in S$, then $s'r=0$ for some element  $s'\in S$. For a left denominator $S$ of the ring $R$, we can form the ring of fractions $S^{-1}R=\{ s^{-1}\, | \, s\in S, \; r\in R\}$.

Suppose that $\OCC$ is {\em a left  denominator set of the associated graded  ring } ${\rm gr}\, R=\bR \oplus \gn / \gn^2\oplus\cdots$ with respect to the $\gn$-adic  filtration. Then the $\OCC$-{\em torsion ideal} of the ring ${\rm gr}\, R$,
\begin{equation}\label{C-tor}
\tau := \tor_{\OCC} ({\rm gr}\, R) = \oplus_{i\geq 1} \tau_i,\;\;
{\rm where}\;\; \tau_i=  \tor_{\OCC} (\gn^i/ \gn^{i+1}),
\end{equation}
  is a homogeneous ideal of the $\N$-graded ring ${\rm gr}\, R$. The factor ring
\begin{equation}\label{C-tor1}
{\rm gr}\, R /\tau = \bR\oplus \gf_1\oplus \gf_2\oplus\cdots ,\;\;
{\rm where}\;\;  \gf_i=(\gn^i/ \gn^{i+1})/\tau_i,
\end{equation}
  is an $\N$-graded ring ($\gf_i\gf_i\subseteq \gf_{i+1}$ for all $i,j\geq 1$) and a subring of the localization ring
  $$\OCC^{-1} {\rm gr}\, R\simeq \OCC^{-1} ({\rm gr}\, R/\tau )=\bQ\oplus \OCC^{-1}\gf_1\oplus \OCC^{-1}\gf_2\oplus\cdots $$
  which is an $\N$-graded ring.

Suppose, in addition, that the nilpotency degree $\nu$ of the prime radical $\gn$ is {\em finite}. Then the prime radical $\gn_{{\rm gr}\, R/\tau }$ of the ring ${\rm gr}\, R/\tau $ is equal to
\begin{equation}\label{NgrR}
\gn_{{\rm gr}\, R/\tau }=\gf := \oplus_{i\geq 1} \gf_i.
\end{equation}
It is a nilpotent ideal of nilpotentcy degree $\max \{ i\geq 1\, | \, \gf_i\neq 0\} \leq \nu $.
\begin{theorem}\label{3Sep12}
Let $R$ be a ring. The following statements are equivalent.
\begin{enumerate}
\item The ring $R$ has a left Artinian ring left quotient ring $Q$.
\item The set $\OCC$ is a left denominator set in the ring ${\rm gr}\, R$, $\OCC^{-1} {\rm gr}\, R$ is a left Artinian ring,  $\gn$ is a nilpotent ideal and $\CC'\subseteq \CC$.
\item The set $\OCC$ is a left denominator set in the ring ${\rm gr}\, R$,  the left quotient ring $Q({\rm gr}\, R/ \tau )$ of the ring ${\rm gr}\, R/\tau$ is a left Artinian ring, $\gn$ is a nilpotent ideal and $\CC'\subseteq \CC$.
\end{enumerate}
If one of the equivalent conditions holds then ${\rm gr}\, Q\simeq Q({\rm gr}\, R/ \tau )\simeq \OCC^{-1} {\rm gr}\, R$ where ${\rm gr}\, Q$ is the associated graded ring with respect to the prime radical filtration.
\end{theorem}

{\it Proof}. $(1\Rightarrow 2)$ Suppose that the ring $Q$ is a left Artinian ring. Then so is the associated graded ring ${\rm gr}\, Q$.
 By Theorem \ref{2Sep12}, the ideal $\gn$ is nilpotent, and $\CC'\subseteq \CC$. Let us show that {\em the set} $\OCC$ {\em is a left Ore set in the ring} ${\rm gr} \, R$.
  It suffices to show that, for elements $\bc\in \OCC$ and $ \alpha = a+\gn^i \in \gn^i / \gn^{i+1}$, there exist elements $\bc_1\in \OCC$ and $\alpha_1= a_1+\gn^i\in \gn^i / \gn^{i+1}$ such that
$$ \bc_1 \alpha = \alpha_1 \bc.$$
If $\alpha \in \tau_i$ then $\bc_1\alpha =0$ and it suffices to take $\alpha_1=0$. If $\alpha \not\in \tau_i$ then $\bs \alpha\neq 0$ for all elements $ \bs \in \OCC$. Recall that the set $\CC'= \CC$ is a left Ore set in $R$ (Theorem \ref{2Sep12}) and $c\in \CC$ (by the condition (c) of Theorem \ref{2Sep12}) where $\bc = c+\gn$. So, there are elements $c_1\in \CC'$ and $a_1\in R$ such that
$$ c_1a=a_1c.$$
Fix the elements $c_1\in \CC'$ and $a_1\in \gn^j\backslash \gn^{j+1}$ with maximal possible degree $j$ of the element $a_1$ (clearly, $a_1\neq 0$ and $j\leq i$ since $\bc_1\alpha \neq 0$). We claim that  $j=i$. Suppose not, i.e. $j<i$, we seek a contradiction. Then
$$ 0\equiv c_1a\equiv a_1c\mod \gn^{j+1}.$$ By the condition (e) of Theorem \ref{2Sep12}, there exists an element $c'$ such that $c'a_1\equiv 0\mod \gn^{j+1}$, i.e. the degree of the element $c'a_1$ is strictly larger than $j$ but the equality holds $$ c'c_1a= c'a_1c,$$
  which contradicts the maximality of $j$.  Therefore, $j=i$. It suffices  to take the equality $c_1a=a_1c\in \gn^i$ modulo $\gn^{i+1}$ to obtain the result, $\bc_1\alpha =\alpha_1\bc\in \gn^i \in \gn^i/\gn^{i+1}$.

By the condition (e) of Theorem \ref{2Sep12}, the set $\OCC$ is a left denominator set in the ring ${\rm gr} \, R$.
 To finish the proof, it suffices to show that ${\rm gr}\, Q\simeq \OCC^{-1}{\rm gr}\, R$ (since ${\rm gr}\, Q$ is a left Artinian ring). In brief, this follows from statements (A2), (A3) and (C) established in the proof of Theorem \ref{2Sep12}. We keep the notation of Section \ref{NCLALQ}. Recall that

(C) $\CC = \CC'$ and $Q=Q'$.

Then statements (A2) and (A3) can be written as follows

(A2') $(\CC^{-1}\gn )^i=\CC^{-1}\gn^i$  {\em for all } $i\geq 1$;

(A3') {\em the map}
$$\CC^{-1}R/\CC^{-1}\gn\ra \bQ, \;\; c^{-1}r\mapsto \bc^{-1} \br, $$ {\em is a ring isomorphism with the inverse} $\bc^{-1} \br\mapsto c^{-1}r$.

By (A2'), the ideal $\CC^{-1}\gn$ of the ring $Q$ is a nilpotent ideal such that the factor ring
$$ Q/\CC^{-1}\gn \simeq \bQ$$ is a semisimple  (by (A3')), hence semiprime. Therefore, the ideal $\CC^{-1}\gn$ is {\em the prime radical of the ring} $Q$. Now,
\begin{eqnarray*}
 {\rm gr} \, Q& = & Q/\CC^{-1}\gn \oplus\cdots  \oplus (\CC^{-1}\gn )^i/(\CC^{-1}\gn )^{i+1}\oplus\cdots \stackrel{(A2'), (A3')} {=}
  \bQ \oplus\cdots  \oplus \CC^{-1}\gn ^i/\CC^{-1}\gn ^{i+1}\oplus\cdots\\
 &\simeq &  \bQ \oplus\cdots  \oplus \CC^{-1}(\gn ^i/\gn ^{i+1})\oplus\cdots \simeq \bQ \oplus\cdots  \oplus \OCC^{-1}(\gn ^i/\gn ^{i+1})\oplus\cdots\\
 &\simeq & \OCC^{-1}{\rm gr}\, R \simeq \OCC^{-1}({\rm gr}\, R/\tau ).
\end{eqnarray*}

Since the ring ${\rm gr}\, Q$ is a left Artinian ring (since $Q$ is so), every regular element of ${\rm gr}\, Q$  is a unit. This is also true for the ring $\OCC^{-1}{\rm gr}\, R\simeq \OCC^{-1} ({\rm gr}\, R/\tau )$ (since ${\rm gr}\, Q\simeq \OCC^{-1}{\rm gr}\, R$). Therefore, $\OCC^{-1}{\rm gr}\, R= Q({\rm gr}\, R/\tau )$ is  the left quotient ring of the ring ${\rm gr}\, R/\tau $.

$(2\Rightarrow 3)$ $Q({\rm gr}\, R/\tau )=Q(\OCC^{-1}{\rm gr}\, R)=\OCC^{-1}{\rm gr}\, R$ since the ring $\OCC^{-1}{\rm gr}\, R$ is a left Artinian ring (any left Artinian ring coincides with its left quotient ring).

$(3\Rightarrow 1)$ It suffice to show that the conditions (a)--(e) of Theorem \ref{2Sep12} are satisfied. The conditions (b) and (c) are given. By (\ref{NgrR}), the prime radical $\gn_{{\rm gr}\, R/ \tau}$ of the ring ${\rm gr}\, R/ \tau$ is equal to $\gf = \oplus_{i\geq 1} \gf_i$ and $({\rm gr}\, R/ \tau)/\gn_{{\rm gr}\, R/ \tau}\simeq     \bR$. By the assumption,  the ring $Q({\rm gr}\, R/ \tau )$ is a left Artinian ring, and so the conditions (a)--(e) of Theorem \ref{2Sep12} are satisfied for the ring ${\rm gr}\, R/ \tau $. Since the set $\OCC$ is a left Ore set of the ring ${\rm gr}\, R/ \tau$ and  consists
 of regular elements, the conditions (a), (d) and (e) for the ring ${\rm gr}\, R/ \tau$ are the same as the conditions (a), (d) and (e) for the ring $R$. $\Box $



\section{Criteria similar to  Robson's Criterion}\label{WVRC}

In this section, two criteria similar to  Robson's Criterion are found (Theorem \ref{4Sep12} and Corollary \ref{a10Sep12}):  Robson's Criterion holds where $\CC$ is replaced by $\CC'$ and one of the conditions is changed accordingly  (Theorem \ref{4Sep12}), Corollary \ref{a10Sep12} is a `weaker' version of Theorem \ref{4Sep12}.

{\bf Robson's Criterion}. A ring $R$ is called $\gn$-{\em reflective} if, for $c\in R$, $c$ is regular in $R$ iff $c+\gn$ is regular in $\bR$; equivalently, $\CC'=\CC$. A ring $R$ is called $\gn$-{\em quorite} if, given $c\in \CC$ and $n\in \gn$, there exist $c'\in \CC$ and $n'\in \gn$ such that $c'n=n'c$.
 A left ideal $I$ of the ring $R$ is called a $\CC$-{\em closed} if, for elements $c\in \CC$ and $r\in R$, $cr\in I$ implies $r\in I$. Similarly, a $\CC'$-{\em closed} right  ideal is defined.

\begin{theorem}\label{Th-2.10-Robson}
{\rm (Robson, Theorem 2.10, \cite{Robson-ArtQuotRings-67})} Let $R$ be a ring. The following statements are equivalent.
\begin{enumerate}
\item The  ring $R$ has  a left Artinian left quotient ring $Q$.
\item The ring $R$ is $\gn$-reflective and $\gn$-quorite, $\bR$ is a left Goldie ring,  $\gn$ is a nilpotent ideal and the ring $R$ satisfies ACC on $\CC$-closed left ideals.
\end{enumerate}
\end{theorem}

The next result shows that the condition `$R$ is $\gn$-reflective' can be weaken.

\begin{theorem}\label{4Sep12}
Let $R$ be a ring. The following statements are equivalent.
\begin{enumerate}
\item The  ring $R$ has  a left Artinian left quotient ring $Q$.
\item
\begin{enumerate}
\item The ring $\bR $ is a left Goldie ring.
\item $\gn$ is a nilpotent ideal.
\item $\CC'\subseteq \CC$.
\item  If $c\in \CC'$ and $n\in \gn$ then there exist elements $c_1\in \CC'$ and $n_1\in \gn$ such that $c_1n=n_1c$.
\item The ring $R$ satisfies ACC for  $\CC'$-closed left  ideals.
\end{enumerate}
\end{enumerate}
\end{theorem}

{\it Proof}. $(1\Rightarrow 2)$ If statement 1 holds then $\CC = \CC'$ (Theorem \ref{2Sep12}) and the conditions (a)--(e) are precisely Robson's Criterion (Theorem 2.10, \cite{Robson-ArtQuotRings-67}).

$(2\Rightarrow 1)$  The proof of this implication has a similar structure to the proof of the implication $(2\Rightarrow 1)$ of Theorem \ref{2Sep12}. It comprises three steps:

(A) $\CC'$ {\em is a left Ore set in} $R$. From which it follows easily that

(B) {\em the ring $Q':= \CC'^{-1}R$ is a left Artinian ring}, and then that

(C) $\CC= \CC'$ {\em and} $Q=Q'$.

(A) $\CC'$ {\em is a left Ore set in} $R$. We have to show that, for elements $c\in \CC'$ and $r\in R$, there exist elements $c'\in \CC'$ and $r'\in R$ such that
$$ c'r=r'c.$$
To prove this statement we use a downward induction on the  degree $i= \deg (r)$ of the element $r$. The base of the induction $i=\nu +1$, i.e. $r=0$, is obvious (take $c'=1$ and $ r'=0$). Suppose that $i<\nu +1$ and that statement (A) is true for all elements of the ring $R$ with degrees $i'$ such that $i<i'\leq \nu +1$.

{\em Case 1:} $i=0$. The ring $\bR$ is a semiprime left Goldie ring. By Goldie's Theorem,
$$ \bc_1\br =\br_1\bc$$ for some elements $\bc_1= c_1+\gn \in \OCC$ and $\br_1=r_1+\gn\in \bR$. By the condition (c), $c_1\in \CC'$. Then, the element
$$ n:= c_1r-r_1c\in \gn$$
has degree $\geq 1$, and, by induction, there are elements $c_2\in \CC'$ and $r_2\in R$ such that $c_2n=r_2c$, and so
$$ c_2c_1r=(r_2+c_2r_1)c.$$
It suffices to take $c'=c_2c_1$ and $r'=r_2+c_2r_1$.

{\em Case 2:} $i>0$. Applying the condition (d),  we see that if $c\in \CC'$  and $r\in \gn^i$ then there exist elements $c'\in \CC'$ and $r'\in \gn$ such that
\begin{equation}\label{cn=cn}
c'r=r'c
\end{equation}
The proof of  statement (A) is complete.

 Now, the statements (A1)--(A4) (see the proof of Theorem \ref{2Sep12}) are true  with the same proofs as in the proof of Theorem \ref{2Sep12}.

 (A1) $\gn\CC'^{-1}\subseteq \CC'^{-1}\gn$. In particular, $\CC'^{-1}\gn$ is an {\em ideal}  of the ring $Q'$ such that

(A2) $(\CC'^{-1}\gn )^i=\CC'^{-1}\gn^i$  {\em for all } $i\geq 1$ (obvious), and

(A3) {\em the map}
$$\CC'^{-1}R/\CC'^{-1}\gn\ra \bQ, \;\; c^{-1}r\mapsto \bc^{-1} \br, $$ {\em is a ring isomorphism with the inverse} $\bc^{-1} \br\mapsto c^{-1}r$.

{\it Proof of (A3)}: By (A1), the condition (c) and the universal property of left Ore localization, the both
 ring homomorphisms in (A3) are well-defined; clearly they are mutually inverse.

(A4) $\CC'^{-1}\gn$ {\em is the prime radical of the ring} $Q'$.

 (B) {\em The ring $Q'$ is a left Artinian ring}. By (A3) and Goldie's Theorem, the ring $$\CC'^{-1}R/\CC'^{-1}\gn \simeq \bQ$$ is a {\em semisimple ring}.
 For all integers $i\geq 1$, the left $\bQ$-modules/$Q'$-modules $\CC'^{-1}\gn^i/\CC'^{-1}\gn^{i+1}\simeq (\CC'^{-1}\gn )^i/(\CC'^{-1}\gn )^{i+1}$ are semisimple, hence of finite length, by the condition (e) (by (A), the condition (e) is equivalent to the condition that the ring $\CC'^{-1}R$ is left Noetherian). Hence, the ring $Q'$ is a left Artinian ring. Then

 (C) $\CC=\CC'$ {\em and} $Q=Q'$ (repeat the proof of statement (C) in the proof of Theorem \ref{2Sep12}).  $\Box$

$\noindent $

The next corollary shows that the condition (c) in Theorem \ref{4Sep12} can be weaken.

\begin{corollary}\label{a10Sep12}
Let $R$ be a ring. The following statements are equivalent.
\begin{enumerate}
\item The  ring $R$ has  a left Artinian left quotient ring $Q$.
\item
\begin{enumerate}
\item The ring $\bR$ is a left Goldie ring.
\item $\gn$ is a nilpotent ideal.
\item  There exists a submonoid $\CC''$ of $\CC$ such that
$\overline{\CC''}= \CC$.
\item  If $c\in \CC''$ and $n\in \gn$ then there exist elements $c_1\in \CC''$ and $n_1\in \gn$ such that $c_1n=n_1c$.
\item The ring $R$ satisfies ACC for  $\CC''$-closed left  ideals.
\end{enumerate}
\end{enumerate}
If one of the equivalent conditions holds then $\CC''$ is a left Ore  set in $R$, $\CC''^{-1}R=Q$, $\CC''^{-1}\gn$ is the prime radical of the ring $Q$ which is a nilpotent ideal of nilpotency degree $\nu$, and the map $Q/ \CC''^{-1}\gn \ra \bQ$, $c^{-1}r\mapsto \bc^{-1}\br$,  is a ring isomorphism with the inverse $\bc^{-1}\br\mapsto c^{-1}r$ where $c$ is any element of $\CC''$ such that $\bc = c+\gn$.
\end{corollary}

{\it Proof}. $(1\Rightarrow 2)$ Take  $\CC = \CC'$ and use Theorem \ref{4Sep12}.

$(2\Rightarrow 1)$  The proof of this implication proceeds in a similar manner as the proof of the implication $(2\Rightarrow 1)$ of Theorem \ref{4Sep12}.

(A'') $\CC''$ {\em is a left Ore set in} $R$. From which it follows easily that

(A1'') $\gn\CC''^{-1}\subseteq \CC''^{-1}\gn$.

(A2'') $(\CC''^{-1}\gn )^i=\CC''^{-1}\gn^i$  {\em for all } $i\geq 1$.

(A3'') {\em The map}
$$\CC''^{-1}R/\CC''^{-1}\gn\ra \bQ, \;\; c^{-1}r\mapsto \bc^{-1} \br, $$ {\em is a ring isomorphism with the inverse} $\bc^{-1} \br\mapsto c^{-1}r$ {\em where $c$ is any element of $\CC''$ such that} $\bc = c+\gn$.

(A4'') $\CC''^{-1}\gn$ {\em is the prime radical of the ring} $Q''$.

The proofs of these statements are identical (apart from obvious minor modifications) to those of statements (A) and (A1) -- (A4) in the proof of Theorem \ref{4Sep12}.

(B'') {\em the ring $Q'':= \CC''^{-1}R$ is a left Artinian ring}. By (A3'') and Goldie's Theorem, the ring $\CC''^{-1}R/\CC''^{-1}\gn\simeq \bQ$ is a semisimple ring. For all integers $i\geq 1$, the left $\bQ$-modules/$Q''$-modules $\CC''^{-1}\gn^i/\CC''^{-1}\gn^{i+1}\simeq (\CC''^{-1}\gn )^i/(\CC''^{-1}\gn )^{i+1}$ are semisimple, hence of finite length, by the condition (e) (by (A''), the condition (e) is equivalent to the condition that the ring $\CC''^{-1}R$ is left Noetherian). Therefore, the ring $Q''$ is a left Artinian ring.

Proofs of the statements below are identical to the same statements in the proof of Theorem \ref{C2Sep12}.

(B1'') $\CC\subseteq Q''^*$.

(B2'') $\CC$ {\em is a left Ore set in} $R$.


 (C'') $Q=Q''$.  $\Box$


\section{A left quotient ring of a factor ring}\label{LQFAR}

The aim of this section is to prove Theorem \ref{13Sep12} which, for a ring $R$ with a left Artinian left quotient ring $Q$ and its $\CC$-closed ideal $I\subseteq \gn$, shows that the factor ring $R/I$ has a left Artinian left quotient ring $Q(R/I)$.

\begin{theorem}\label{13Sep12}
Let $R$ be a ring with a left Artinian left quotient ring $Q$, and $I$ be a $\CC$-closed ideal of $R$ such that $I\subseteq \gn$. Then
\begin{enumerate}
\item The set $\CC_{R/I}$ of regular elements of the ring $R/I$ is equal to the set $\{ c+I\, | \, c\in \CC \}$.
\item The ring $R/I$ has a left Artinian left quotient ring $Q(R/I)$, $\CC^{-1}I$ is an ideal of $Q$  and the map $Q/\CC^{-1}I\ra Q(R/I)$, $ c^{-1}r+\CC^{-1}I\mapsto (c+I)^{-1}(r+I)$, is a ring isomorphism with the inverse $(c+I)^{-1}(r+I)\mapsto c^{-1}r+\CC^{-1}I$.
\end{enumerate}
\end{theorem}

Before giving a proof of Theorem \ref{13Sep12} we need some results that are used in the proof (and also in the proof of Corollary \ref{a14Sep12}).  The next result describes the set of regular elements of a ring and the group of units of its left quotient ring.
\begin{lemma}\label{a13Sep12}
Let $R$ be a ring such that the set of regular elements $\CC$ is a left Ore set, $Q^*$ be the group of units of the left quotient ring $Q= \CC^{-1}R$ of the ring $R$. Then
\begin{enumerate}
\item $C= R\cap Q^*$.
\item $Q^*=\{ s^{-1}t\, | \, s,t \in \CC\}$.

\end{enumerate}
\end{lemma}

{\it Proof}. 1. Clearly, $R\cap Q^*\subseteq \CC$ and $C\subseteq R\cap Q^*$ (since $C\subseteq Q^*$), i.e. $\CC = R\cap Q^*$.

2. Clearly, $Q^* \supseteq \{ s^{-1}t\, | \, s,t\in \CC\}$. Let $q\in Q^*$. Then $q=s^{-1}t$ for some elements $s\in \CC$ and $t\in R$. Then $t=sq\in R\cap Q^* = C$, by statement 1. Therefore, $Q^*=\{ s^{-1}t\, | \, s,t,\in \CC\}$. $\Box $

$\noindent $

A left/right ideal of a ring is called {\em nil} if it consists of nilpotent elements. A nilpotent left/right ideal is nil but not vice versa, in general.

\begin{lemma}\label{a17Jul12}
Let $I$ be a nil ideal of a ring $R$;  $\xi :R\ra R/I$, $r\mapsto r+I$; $R^*$ and $(R/I)^*$ be the groups of units of the rings  $R$ and $R/I$ respectively. Then
\begin{enumerate}
\item The set $1+I$ is a normal subgroup of the group $R^*$.
\item The short exact sequence of group homomorphisms $1\ra 1+I\ra R^*\stackrel{\xi^* }{\ra} (R/I)^*\ra 1$ is exact where $\xi^*=\xi|_{R^*}$.
\item Let $\bu=u+I\in (R/I)^*$ where $u\in R$. Then $(\xi^*)^{-1} (\bu ) = u+I= \xi^{-1}(u)$. In particular, $R^* +I = R^*$ and $\xi^{-1}((R/I)^*)=R^*$.

\end{enumerate}
\end{lemma}

{\it Proof}. 1. If $a\in I$ then $a^n=0$ for some natural number $n$, and so $(1-a)^{-1} = 1+\sum_{i=1}^{n-1}a^i\in 1+I$. Therefore, $1+I$ is a subgroup of the group $R^*$. For all elements $u\in R^*$ and $a\in I$, $u(1+a)u^{-1} = 1+uau^{-1}\in 1+I$, i.e. the subgroup $1+I$ is a normal subgroup of $R^*$.

2 and 3.  The map $\xi^*$ is a group homomorphism with $\ker (\xi^*) = 1+I$, by statement 1. Let us show that  the map $\xi^*$ is a surjection. Let $\bu=u+I\in (R/I)^*$ where $u\in R$. We have to show that $u\in R^*$. Notice that $\bu^{-1} = v+I$ for some element $v\in R$. Then $uv=1+a$ and $vu=1+b$ for some elements $a,b\in I$. By statement 1, the elements $1+a$ and $1+b$ are units of the ring $R$. Therefore, $u\in R^*$, and so $(\xi^*)^{-1} (\bu ) = u+I$. $\Box $


\begin{proposition}\label{a20Jul12}
Let $R$ be a ring, $I$ be an ideal of the ring $R$, and $S$ be a left Ore set in $R$ that consists of regular elements of the ring $R$. If $S^{-1}R$ is a left Artinian (resp. left Noetherian) ring then $S^{-1}I$ is an ideal of the ring $S^{-1}R$.
\end{proposition}

{\it Proof}. We have to show that that the left ideal $S^{-1}I$ of the ring $S^{-1}R$ is also a right ideal, i.e.
$$S^{-1}Is^{-1} \subseteq S^{-1}I\;\; {\rm  for\;\; all}\;\;  s\in S.$$ Since the ring $S^{-1}R$ is a left Artinian (resp.  left Noetherian) ring the chain of left ideals in $S^{-1}R$: $S^{-1}Is\supseteq S^{-1}Is^2\supseteq \cdots \supseteq S^{-1}Is^n\supseteq \cdots $ (resp. $S^{-1}Is^{-1}\subseteq S^{-1}Is^{-2}\subseteq \cdots \subseteq S^{-1}Is^{-n}\subseteq \cdots $)  stabilizers say on $n$'th step,  $S^{-1}Is^n= S^{-1}Is^{n+1}$ (resp. $S^{-1}Is^{-n}= S^{-1}Is^{-n-1}$). In both cases, we have the equality $S^{-1}Is^{-1}= S^{-1}I$ since $s$ is a regular element of the ring $R$.  $\Box $

$\noindent $

Let $R$ be a ring with a left Artinian left quotient ring $Q$. Then, the left $Q$-module $Q$ has finite length, its (Jacobson) radical $\gr$ is a nilpotent ideal of nilpotency degree, say,  $\nu_Q$ and coincides with the prime radical of the ring $Q$. Let $\Max (Q) = \{ \gm_1, \ldots , \gm_s\}$ be the set of maximal ideals of the ring $Q$ and $\gp_1 := R\cap \gm_1, \ldots , \gp_s:= R\cap \gm_s$. Then $\gr = \cap_{i=1}^s\gm_i$ since $\Min (Q)= \{ \gm_1, \ldots , \gm_s\}$ is  the set of minimal primes   of the ring $Q$.

\begin{proposition}\label{b13Sep12}
Let $R$ be a ring with a left Artinian left quotient ring $Q$. Then
\begin{enumerate}
\item If $I$ is an ideal of the ring $R$ then $\CC^{-1}I$ is an ideal of the ring $Q$, $(\CC^{-1}I)^i= \CC^{-1}I^i$ for all $i\geq 1$.
\item $\Min (R) = \{ \gp_1, \ldots , \gp_s\}$, the set of minimal  prime ideals of the ring $R$; and $\gn =\cap_{i=1}^s\gp_i$.
\item $\gn = R\cap \gr$, $\CC^{-1} \gn =\gr$, $\gn$ is a nilpotent ideal of $R$ of nilpotency  degree $\nu =\nu_Q$.
    \item $\CC' = \CC$, i.e. if $c\in R$ then $c\in \CC$ iff $c+\gn \in \OCC$.
            \item The ring $\bR$ is a semiprime left Goldie ring and  the set $\OCC$ of regular elements of $\bR$ is equal to $\{ c+\gn \, | \, c\in \CC \}$.
                \item The map $Q/ \gr \ra \bQ$, $c^{-1}r+\gr \mapsto \bc^{-1}\br$, is a ring isomorphism with the inverse $\bc^{-1}\br \mapsto c^{-1}r+\gr$.
\end{enumerate}
\end{proposition}

{\it Proof}. 1. By Proposition \ref{a20Jul12}, $\CC^{-1}I$ is an ideal of the ring $Q$. In particular, $I\CC^{-1}\subseteq \CC^{-1}I$. Then $(\CC^{-1}I)^i= \CC^{-1}I^i$ for all $i\geq 1$.

2. For every $i=1, \ldots , s$,
$$\CC^{-1}\gp_i = \CC^{-1}(R\cap \gm_i)=\CC^{-1}R\cap \CC^{-1}\gm_i= Q\cap \gm_i = \gm_i.$$
  Therefore, the ideals $\gp_i$ are distinct and $\gp_i\not\subseteq \gp_j$ for all $i\neq j$. The ideals $\gp_i$ are prime: if $\ga\gb\subseteq \gp_i$ for some ideals $\ga$ and $\gb$ of the ring $R$ then $\CC^{-1}\ga$ and $\CC^{-1}\gb$ are ideals of the ring $Q$ (statement 1) such that $\CC^{-1}\ga \CC^{-1}\gb \subseteq \CC^{-1} \ga \gb \subseteq \CC^{-1}\gp_i = \gm_i$; hence either $\CC^{-1}\ga \subseteq \gm_i$ or $\CC^{-1}\gb \subseteq \gm_i$, and so either $\ga \subseteq R\cap \CC^{-1}\ga \subseteq R\cap \gm_i= \gp_i$ or $\gb \subseteq R\cap \CC^{-1}\gb \subseteq R\cap \gm_i= \gp_i$. The intersection $$\cap_{i=1}^s\gp_i = \cap_{i=1}^s (R\cap \gm_i) = R\cap \cap_{i=1}^s \gm_i = R\cap \gr \subseteq \gr$$ is a nilpotent ideal of the ring $R$. Therefore, $\Min (R) = \{ \gp_i , \ldots , \gp_s\}$(by statement 1 and $\Min (Q) = \{ \gm_1, \ldots , \gm_s\}$)  and $\gn = \cap_{i=1}^s \gp_i$.

3. $\gn = \cap_{i=1}^s \gp_i = R\cap \gr$ (see the proof of statement 2).
$$\CC^{-1}\gn = \CC^{-1} (R\cap \gr ) = \CC^{-1} R\cap \CC^{-1} \gr = Q\cap \gr = \gr .$$ By statement 1, for all $i\geq 1$, $\CC^{-1} \gn^i = (\CC^{-1} \gn )^i= \gr^i$. Therefore, $\gn$ is a nilpotent ideal of the ring $R$ of nilpotency degree $\nu_Q$.

4--6. Let us show that $\CC + \gn \subseteq \CC$. Let $c\in \CC$  and $n\in \gn$. We have to show that $c+n\in \CC$. Clearly,  $c^{-1} n \in \CC^{-1} \gn = \gr$ (statement 3), $c^{-1}n$ is a nilpotent element of $Q$ and $1+c^{-1}n\in Q^*$. Then
$$c+n = c(1+c^{-1} n)\in Q^*, $$ hence $c+n\in R\cap Q^* = \CC$, by Lemma \ref{a13Sep12}.(1), as required.

Let, for moment, $\CC_{\bR}$ be the set of regular elements of the ring $\bR$ and $\OCC :=\{ c+n \, | \, c\in \CC\}$. Notice that statement 4 is equivalent  to two statements: $\CC +\gn \subseteq \CC$ and $\CC_{\bR}=\OCC$. Since $\bR = R/ \gn = R/R\cap \gr \subseteq Q/ \gr$ and $\CC\subseteq Q^*$, we see that $ \OCC\subseteq (Q/ \gr)^*$. Therefore,  $\OCC \subseteq \CC_{\bR}$ and $\OCC$ is a left Ore set of the ring $\bR$ that consists of regular elements of $\bR$ (since $\CC$ is a left Ore set of $R$). Notice that
$$ Q/ \gr = \CC^{-1}R / \CC^{-1} \gn \simeq \CC^{-1} (R/ \gn ) \simeq \OCC^{-1} \bR.$$ So, the map $Q/ \gr \ra \OCC^{-1} \bR$, $ c^{-1} r +\gr  \mapsto \bc^{-1} \br$, is a ring isomorphism (where $c\in \CC$ and $ r\in R$). Then $\CC_{\bR}\subseteq (Q/ \gr )^*$ (for every element $\alpha \in \CC_{\bR}$, the map $\cdot \alpha : \OCC^{-1}\bR\ra \OCC^{-1} \bR$, $u\mapsto u \alpha$, is a left $\OCC^{-1}\bR$-module monomorphism, hence an isomorphism since the left $\OCC^{-1}\bR$-module $\OCC^{-1}\bR\simeq Q/\gr$ has finite length; the map $(\cdot \alpha )^{-1}$ has necessarily the form $\cdot q$ for some element $q\in Q/\gr $ such that $\alpha q = q\alpha =1$, i.e. $\alpha \in (Q/\gr )^*$). It follows that $\CC_{\bR}$ is a left Ore set of the ring $\bR$ such that
$$\CC_{\bR}^{-1} \bR = \OCC^{-1} \bR = Q/\gr$$
 (for $\alpha \in \CC_{\bR}$ and $ \br \in \bR$, $\alpha^{-1}\in (Q/ \gr)^*$, and so $\br \alpha^{-1} = \bc^{-1}\br_1$ for some elements $\bc \in \OCC\subseteq \CC_{\bR}$ and $ r_1\in R$, i.e. $\CC_{\bR}$ is a left Ore set of $\bR$). Now, statement 6 is obvious; and $\bR$ is a semiprime left Goldie ring (by Goldie's Theorem).

It remains to show that $\CC_{\bR}\subseteq \OCC$. Let $\alpha \in \CC_{\bR}$. Then $\alpha \in (Q/ \gr)^*$ and, by Lemma \ref{a17Jul12}.(2), $\alpha = q+\gr$ for some $q\in Q^*= \{ s^{-1}t\, | \, s,t \in \CC\}$ (Lemma \ref{a13Sep12}.(2)). Then $q=s^{-1}t$ for some elements $s,t\in \CC$.  Let $\alpha = c+\gn$ for some $c\in R$. It suffices to show that $c\in \CC$. Clearly,
$$n:= sc-t\in R\cap \gr = \gn $$ (Proposition \ref{b13Sep12}.(3)). Notice that $t+n\subseteq \CC+\gn \subseteq \CC$, and then  $ c= s^{-1} (t+n) \in R\cap Q^* = \CC$ (Lemma \ref{a13Sep12}. (1)). $\Box $

$\noindent $

{\bf  Proof of Theorem \ref{13Sep12}}. By Proposition \ref{a20Jul12}, the map
$$\v : Q\ra Q/\CC^{-1}I, \;\;  q\mapsto q+\CC^{-1}I,$$ is a ring epimorphism. The ideal $\CC^{-1}I$ of the ring $Q$ is nilpotent ideal (since $\CC^{-1}I\subseteq \CC^{-1}  \gn = \gr $).  By Lemma \ref{a17Jul12}.(2), the map $\v$ induces the group epimorphism
$$\v^*= \v |_{Q^*}: Q^*\ra (Q/\CC^{-1}I)^*.$$  The ideal $I$ is $\CC$-closed,  i.e. $I=R\cap \CC^{-1}I$. Since the ring $Q/ \CC^{-1}I$ is a left Artinian ring and $R/I= R/R\cap \CC^{-1}I \subseteq Q/ \CC^{-1}I$, we
 see that
 $$C_{R/I}\subseteq (Q/ \CC^{-1}I)^*.$$ The inclusion $C\subseteq Q^*$ implies the inclusion $\v^* (\CC) \subseteq (Q/ \CC^{-1}I)^*$. Therefore,
 $$ \widetilde{\CC}:=\{ c+I\, | \, c\in \CC\} \subseteq R/ I \cap (Q/ \CC^{-1}I)^*\subseteq \CC_{R/I}.$$
 We claim that the equality holds,
\begin{equation}\label{CRIC}
\CC_{R/I}= \widetilde{\CC}.
\end{equation}
It suffices to show that if $\alpha = c+I\in \CC_{R/I}$ (where $c\in R$) then $ c\in \CC$. Using the facts that $\CC_{R/I}\subseteq (Q/\CC^{-1}I)^* = \v^* (Q)$ and $Q=\{ s^{-1}t \, | \, s,t\in \CC\}$ (Lemma \ref{a13Sep12}.(1)), we see that $\alpha = (s+\CC^{-1}I)^{-1} (t+\CC^{-1}I)$ for some elements $s,t\in \CC$; equivalently, there exists  an element $\l \in \CC$ such that
$$ i:= \l ( sc-t)\in I.$$
 Notice that $\l t +i\in \CC +I \subseteq \CC +\gn \subseteq \CC$ (by Proposition \ref{b13Sep12}.(4)). Now,
 $$ c= (\l s)^{-1} (\l t +i) \in R\cap Q^*= \CC,$$
 by Lemma \ref{a13Sep12}.(1), as required. The proof of statement 1 is complete.

 The set $\CC$ is a left Ore set in $R$, hence its image in the ring $Q/\CC^{-1} I$, which is $\widetilde{\CC}= \CC_{R/I}$ (by (\ref{CRIC})), is a left Ore set. Therefore, the quotient ring $Q(R/I):= \CC_{R/I}^{-1} (R/I)$ exists and, by the universal property of left localization,  the map
 $$ \th : Q(R/I) \ra Q/ \CC^{-1} I , \;\; (c+I)^{-1} (r+I)\mapsto (c+\CC^{-1} I)^{-1} (r+\CC^{-1}I) = c^{-1}R+\CC^{-1} I, $$
 is a ring monomorphism (since $R/I\subseteq Q/ \CC^{-1}I$ and $ \widetilde{\CC}= \CC_{R/I} \subseteq (Q/ \CC^{-1}I)^*$). By (\ref{CRIC}), the map $\th$ is a surjection. Therefore, $\th$ is an isomorphism.  In particular, the ring $Q(R/I)$ is a left Artinian ring. The proof of statement 2 is complete.  $\Box $


\begin{corollary}\label{a14Sep12}
The implication $(1\Rightarrow 2)$ of Theorem \ref{2Sep12} holds.
\end{corollary}

{\it Proof}. The conditions (a)--(c)  of Theorem \ref{2Sep12} are the conditions (3)--(5) of Proposition \ref{b13Sep12}. The ring $Q$ is a left Artinian ring. Hence, the left  $Q$-module/$\bQ$-module $\CC^{-1} \gn^i/ \CC^{-1} \gn^{i+1}= (\CC^{-1} \gn )^i/ (\CC^{-1} \gn )^{i+1}$ is semisimple of finite length (Proposition \ref{b13Sep12}.(1)). The left $\bR$-module $\gf_i$ is $\OCC$-torsionfree. Therefore, $\gf_i\subseteq \OCC^{-1}\gf_i$ and $\udim_{\bR}(\gf_i) = \udim_{\bR}(\OCC^{-1}\gf_i)$. The chain of $\bR$-module isomorphisms (Proposition \ref{b13Sep12}.(4--6)),
$$ \CC^{-1}\gn^i / \CC^{-1}\gn^{i+1} \simeq \CC^{-1} \gn^i / \CC^{-1} \gt_i \simeq \CC^{-1} (\gn^i / \gt_i) =\CC^{-1} \gf_i \simeq \OCC^{-1} \gf_i, $$
shows that $\udim_{\bR} (\OCC^{-1}\gf_i) = \udim_{\bR} (\CC^{-1}\gn^i / \CC^{-1} \gn^{i+1}) = {\rm length}_{\bQ} (\CC^{-1}\gn^i / \CC^{-1} \gn^{i+1})<\infty$, i.e.  the condition (d)  of Theorem \ref{2Sep12} holds.

By Proposition \ref{b13Sep12}.(4), every element of $\OCC$ is  of type $\bc = c+\gn$  where $c\in \CC$. The map $\cdot \bc : \gf_i \ra \gf_i$, $f\mapsto f\bc$, is an injection as the restriction of the bijection
$$ \cdot c:\,  \CC^{-1} \gn^i / \CC^{-1} \gn^{i+1}\ra \CC^{-1}\gn^i / \CC^{-1} \gn^{i+1}, \;\; x\mapsto xc, $$
to $\gf_i\subseteq \CC^{-1} \gn^i / \CC^{-1} \gn^{i+1}$. So, the condition (e)  of Theorem \ref{2Sep12} holds. $\Box $

$\noindent $

$${\bf Acknowledgements}$$

 The work is partly supported by  the Royal Society  and EPSRC.

\small{

Department of Pure Mathematics

University of Sheffield

Hicks Building

Sheffield S3 7RH

UK

email: v.bavula@sheffield.ac.uk}

\end{document}